\documentclass[11pt]{amsart}

\newcommand{\Fin}{F_{\infty}}
\newcommand{\sub}{\subseteq}

\newcommand{\normal}{\triangleleft}
\newcommand{\W}{{\mathfrak W}}
\newcommand{\V}{{\mathfrak V}}
\newcommand{\Ni}{{\mathfrak N}}
\renewcommand{\O}{{\mathfrak O}}
\newcommand{\U}{{\mathfrak U}}
\newcommand{\A}{{\mathfrak A}}
\newcommand{\K}{{\mathfrak K}}
\newcommand{\E}{{\mathfrak E}}
\newcommand{\B}{{\mathfrak B}}
\newcommand{\G}{{\mathfrak G}}
\newcommand{\X}{{\mathfrak X}}
\newcommand{\Y}{{\mathfrak Y}}

\newcommand{\N}{{\mathbb N}}
\newcommand{\Z}{{\mathbb Z}}

\newcommand{\Q}{{\mathbb Q}}

\newcommand{\bb}{{\rm(B) }}
\newcommand{\bbb}{{\rm(B^*) }}

\newcommand{\1}{\{1\}}

\newcommand{\card}[1]{\mathrm{card}\left( #1 \right)}

\newcommand{\var}[1]{\mathrm{var}( #1 )}
\newcommand{\varr}[1]{\mathrm{var}( #1 )}

\newcommand{\Wr}{\,\mathrm{Wr}\,}
\newcommand{\Wrr}{\,\mathrm{wr}\,}

\newtheorem{Theorem}{\sc Theorem}[section]
\newtheorem{Lemma}[Theorem]{\sc Lemma}
\newtheorem{Corollary}[Theorem]{\sc Corollary}  
 
\theoremstyle{definition}

\newtheorem{Example}[Theorem]{\sc Example}
\theoremstyle{remark}
\newtheorem{Remark}[Theorem]{\sc Remark}

\address{Department of Computer Science \\ 
Yerevan State University\\
375025 Yerevan, Armenia}
\email{mikaelian@e-math.ams.org
\vskip0.1mm {\it Internet:} www.mikaelian.cjb.net}

\begin{document}

\noindent
{\small
Ricerche di Matematica\\ 
Univ. Studi Napoli, Naples, 52 (2003), 1--19 
\vskip10mm
}

\title[Not locally soluble $SI^*$-groups]{Infinitely many not locally soluble $SI^*$-groups}
\author{Vahagn H.~Mikaelian}

\dedicatory{To Shogik Mikaelian on her birthday}

\begin{abstract}
The class of those (torsion-free) $SI^*$-groups which are not locally soluble, has the cardinality of the continuum. Moreover, these groups are not only pairwise non-isomorphic, but also they generate pairwise different varieties of groups. Thus, the set of varieties generated by not locally soluble $SI^*$-groups is of the same cardinality as the set of all varieties of groups. It is possible to localize a variety of groups which contains all groups and varieties constructed. The examples constructed here continue the well known example of a not locally soluble $SI^*$-group built by Hall and by Kov\'acs and Neumann.
\end{abstract}

\maketitle

\section{Introduction}

Of all the types of generalized soluble groups $SI^*$-groups, that is, the groups possessing ascending normal series with abelian factors in some sense are the ``nearest'' to the soluble groups~\cite{RobinsonFinitenessConditions}. This and the fact that all locally soluble groups are $SI$-groups, that is, groups with normal (not necessarily well-ordered) system with abelian factors (see Subsection~\ref{Examples_for_some_other_classes_of_generalized} for definitions), explain the importance of the example of an $SI^*$-group which is {\it not locally soluble}, built independently by Hall~\cite{Hall_Frattini_Subgroups} and by Kov\'acs and Neumann~\cite{KovacsNeumannOnGeneralizedSolubles}. 
The group constructed in~\cite{Hall_Frattini_Subgroups,KovacsNeumannOnGeneralizedSolubles} is a suitably chosen subgroup in the cartesian wreath product $\bigl(M(\Q,F)\Wr C\bigr)\Wr C$, where $C$ is an infinite cyclic group, $\Q$ is the (naturally ordered) set of rational numbers, $F$ is any countable field, and $M(\Q,F)$ is the corresponding McLain group. After the Hall-Kov\'acs-Neumann construction there was no other example of mentioned type in the literature since sixties.

The main aim of this paper is to show that not locally soluble $SI^*$-groups are `many'. Namely, we offer a method of construction of {\it continuum sets of $SI^*$-groups, which are not locally soluble} (following ~\cite{KovacsNeumannOnGeneralizedSolubles} we use the term {\it not locally soluble\/} instead of {\it locally insoluble\/} used, say, in~\cite{RobinsonFinitenessConditions}). Moreover, the groups constructed are not only pairwise non-isomorphic, but also {\it they generate pairwise different varieties of groups}. Taking into account the fact that the cardinality of the set of all varieties of groups  also is continuum, we conclude that {\it the set of varieties generated by not locally soluble $SI^*$-groups is of the same cardinality as the set of all varieties of groups}.

Our construction uses the idea of~\cite{Hall_Frattini_Subgroups,KovacsNeumannOnGeneralizedSolubles} concerning the embedding of an insoluble $SI^*$-group $M(\Q,F)$ into a $2$-generator $SI^*$-group: the latter cannot, therefore, be locally soluble because it contains $M(\Q,F)$. Here we do not use the McLain group but replace it by a group $G$ with properties needed for our construction. We also replace the embeddings of~\cite{Hall_Frattini_Subgroups,KovacsNeumannOnGeneralizedSolubles} by the verbal embedding construction built in~\cite{SubnormalEmbeddingTheorems,SubnormalEmbedOfOrderedGr,OnNormalEmbeddingsOfSoluble} for embeddings of countable groups into two-generator groups with  various pre-given properties. 
  
Section~\ref{section2} below presents our verbal embedding 
construction. In order to avoid a repetition of a large part of the mentioned papers~\cite{SubnormalEmbeddingTheorems,SubnormalEmbedOfOrderedGr,OnNormalEmbeddingsOfSoluble}, we present here only the scheme of our argument. The straightforward application of the latter already enables us to construct {\it countably\/} many not locally soluble $SI^*$-groups (Theorem~\ref{Counably_many}).

In Section~\ref{section3} we modify our initial method and build {\it a continuum} of groups of mentioned type and localize a variety that contains all these groups (Theorems~\ref{central}, \ref{many_varieties}, \ref{new_central} and \ref{new_many_varieties}). This modification is based on a construction of Ol'shanskii~\cite{Olshanski_On_problem-of_finite_base} concerning the continuum of subvarieties in the variety $\G_5 \cap \B_{8qr}$, where $\G_5$ is the variety of soluble groups of length $5$ and $\B_{8qr}$  is the Burnside variety of exponent $8qr$ with $q$ and $r$ being different primes.

In Section~\ref{section4}  we consider possibilities to apply our construction to obtain examples of different nature for other classes of generalized soluble groups: $SI$-groups which are not $SI^*$-groups, $SN$-groups which are not $SN^*$-groups, etc. But for some reasons mentioned there such examples are less interesting than those of not locally soluble $SI^*$-groups and thus we restrict ourselves to an outline of possible modifications. Closing the current paper we show how can the groups constructed also be {\it linearly ordered}. 

For information on the theory of generalized soluble and generalized nilpotent groups we refer to the original articles of Plotkin~\cite{PlotkinClassical} and of Kuro\v s and Chernikov~\cite{KuroshChernikovClassical} as well as to the books of Robinson~\cite{RobinsonFinitenessConditions} and of Kuro\v s~\cite{KuroshThoreyOfGroups}.
Information on wreath products can be found in the paper of Neumann~\cite{OnTheStructureOfWreathProducts}, in the books of Kargapolov and Merzlyakov~\cite{KargapolovMerzlyakov} and of Meldrum~\cite{MeldrumWr}.
For general information on varieties of groups we refer to
the book of Hanna Neumann~\cite{HannaNeumann}. 
Information on linearly (or fully) ordered groups can be found in the book of Fuchs~\cite{BookOfFuchs} or in the papers of Levi~\cite{LeviOnOrderedGroups1,LeviOnOrderedGroups2} and of Neumann~\cite{NeumannOnOrderedGroups,EmbeddingsIntoOrdered}.

\vskip13mm
\section{Verbal embeddings and not locally soluble groups}
\label{section2}

\subsection{The verbal embedding construction}
The verbal embedding mentioned in the introduction is given by:

\begin{Theorem}[See Theorem 1 and Theorem 2 in~\cite{OnNormalEmbeddingsOfSoluble}]
\label{From_That_Paper}
Let $G$ be a countable $SI^*$-group and $V\sub \Fin$ be a non-trivial word set. Then there exists a two-generator $SI^*$-group $H=H(G,V)$ with a subnormal subgroup $\tilde G$ which is isomorphic to $G$ and which lies in the verbal subgroup $V(H)$.
\end{Theorem}

Some details and notations of proof of this theorem are relevant for purposes of Sections~\ref{section3}--\ref{section4} of this paper. Nevertheless, as we said in the introduction, to avoid repetition of a large part of~\cite{SubnormalEmbeddingTheorems,SubnormalEmbedOfOrderedGr} (and in particular of ~\cite{OnNormalEmbeddingsOfSoluble}), we restrict ourselves here by an outline of the proof and refer the reader to that paper for details. Of course, one disadvantage of such an approach is that this paper becomes much dependant, in particular, on~\cite{OnNormalEmbeddingsOfSoluble}.

\begin{proof}[Proof scheme for Theorem~\ref{From_That_Paper}]
Consider complete wreath product $G \Wr S$, where $S$ is a free nilpotent group of some class $c=c(V)$ and of rank $r=r(V)$, such that 

\begin{enumerate}

\item the verbal subgroup $V(S)$ is non-trivial,

\item the group $S$ discriminates~\cite{HannaNeumann} the variety $\var S$.

\end{enumerate}

It is easy to choose such a group $S$. It suffices to take for the given $V$ such a positive integer $c(V)$ that the variety $\Ni_{c(V)}$ does not lie in $\V = \var{\Fin/V(\Fin)}$ and to put $S=F_c(\Ni_{c(V)})$, that is $r(V)=c$. The group $S$ discriminates $\var S = \Ni_{c(V)}$   (see~\cite{HannaNeumann}). (The fact that $S$ is a discriminating group will be used later in the Section~\ref{section3}.)

The verbal subgroup $V(S)$ is non-trivial and contains, thus, an element $a$ of infinite order. Let $\chi_g$ be the elements of the base group $G^S$ of $G\Wr S$ defined as:
$$
  \chi_g (s) =  \left\{
                  \begin{array}{l}
                  g,  \,\,\, \hbox{{\rm if}}\,\,\,  s=a^i, \,\,\, i= 0,1,2,...  ,\\
	1,  \,\,\,  \hbox{{\rm if}}\,\,\, s\in S\, \backslash \{  a^i \,\,| \,\,\, i=0,1,2,...    \}.
                  \end{array}
                  \right.
$$
The subgroup $T=T(G,S)=\langle  \chi_g, S \,|\, g\in G  \rangle$ of $G\Wr S$ contains the first copy of $G$ in $G^S$. Moreover, this first copy lies in the verbal subgroup $V(T)$. 

Now let $\{ G_\delta\,|\, \delta\in \Delta  \}$ be an appropriate ascending  normal soluble series of the $SI^*$-group $G$, and let $L=T \cap G^S$. Then the subgroups 
$$
		M_\delta = L \cap \prod_{s\in S} G_\delta, \quad \delta \in \Delta
$$
do form an ascending normal soluble series for $SI^*$-group $L$. This series can be continued to a series of the group $T$ by addition of preimages of the lower central series  of  nilpotent group $T/L$ under natural endomorphism $T\to T/L$. So $T$ is an $SI^*$-group with ascending soluble normal series $\{ M_\delta\,|\, \delta\in \Delta'  \}$, where the index set $\Delta'$ contains $\Delta$ and at most $c$ extra elements (see~\cite{OnNormalEmbeddingsOfSoluble} for details). 

Modifying this construction we get in~\cite{OnNormalEmbeddingsOfSoluble} an embedding of the group $T$ into an appropriate group $D=D(T,C)$. Namely, we consider the complete wreath product $T \Wr C$, where $C = \langle c \rangle$ is an infinite cycle. Set further for each $g\in T$:
$$
		\pi_g(c^i) =\left\{
                  \begin{array}{l}
                  g, \,\,\, \hbox{{\rm if}}\,\,\,  i \geq 0,  \\
                  1, \,\,\, \hbox{{\rm if}}\,\,\,  i<0.
                  \end{array}
                  \right.
$$
Then the first copy of $T$ in the base group $T^{\,C}$ lies in the derived subgroup of the group
$$
	D= \langle \pi_g , c \,|\,\, g\in T \rangle.
$$
Now we ``lift'' the ascending soluble normal series $\{ M_\delta\,|\, \delta\in \Delta'\}$ of $T$ to an ascending soluble normal series $\{ K_\delta\,|\, \delta\in \Delta''\}$ ($\Delta''$ consists of $\Delta'$ and of  one additional element) of $D$ in the following way. We set $Y=D \cap T^{\,C}$ and define:
$$
	K_\delta = \prod_{c^i\in C}M_\delta \cap Y, \quad \delta \in \Delta'.
$$ 
As the last member of the series for the index $\delta\in \Delta'' \backslash \Delta'$ we add the group $D$ itself.

Now consider the complete wreath product $D \Wr Z$, where $Z = \langle z \rangle$ is an infinite cycle. Clearly $\card D = \card T = \card G = \aleph_0$. So let us write $D$  as 
$$
	D = \{ d_0, d_1, \ldots, d_n, \ldots \,|\, n\in \N \}
$$ 
and define the following element $\omega$ in the base group $D^Z$:
$$
\omega (z^i) = \left\{
                  \begin{array}{ll}
                  d_k,  & {\rm if}\,\,\,  i = 2^k, \,\,\, k=0,1,2,...,  \\
                  1,   &  {\rm if}\,\,\,  i \in \Z \backslash \{ 2^k\,\, | \,\,\,  k=0,1,2,...  \}.
                  \end{array}
                  \right.
$$
Define $H = H(G,V) = \langle \omega,z \rangle$.  $H$ contains $D'=[D,D]$ (and, thus, the groups $T$ and $G$) as {\it subnormal} subgroups. Moreover $G\sub V(H)$. And, finally, we make $H$  an $SI^*$-group in the following way. Let $W_\delta$ be the direct power $\overline{\prod}_{z\in Z} K_\delta$, \, $\delta \in \Delta''$. The subgroups $W_\delta$,  $\delta \in \Delta$, form  as ascending normal soluble series in the direct power $W= \overline{\prod}_{z\in Z} D$.
Intersections
$$
	B_\delta=W_\delta \cap H, \quad \delta \in \Delta''
$$
form an ascending soluble normal series for  $H\cap W$. This series can be continued to a series $\{B_\delta \,|\, \delta \in \Delta'''\}$ for the whole group $H$: we simply add one more member, namely, the group $H$ itself.
\end{proof}

\begin{Remark}
Let us notice that if our initial group $G$ is torsion-free, the group $H(G,V)$ we constructed will also be torsion-free. 
\end{Remark}


A consequence of the proof of Theorem~\ref{From_That_Paper} is:

\begin{Corollary}
\label{UNA2}
If the group $G$ belongs to variety $\U$ and if 
$$
\V = \var{\Fin/V(\Fin)}
$$ 
is the variety corresponding to the word set $V$, than the group $H=H(G,V)$ can be chosen to belong to variety 
$$
	\U \cdot \Ni_{c(V)} \cdot \A^2
$$
(where, as above, $\Ni_{c(V)}\not\sub \V$).
\end{Corollary}

\begin{proof}
Since $T$ is some subgroup of $G\Wr S$, it belongs to the variety $\var G \cdot \var S \sub \U \cdot \Ni_{c(V)}$. The group $D$ is chosen in the variety $\U \cdot \Ni_{c(V)} \cdot \A$. Finally, the group $H$ is a subgroup of wreath product $D \Wr Z$, which belongs to $\var D \cdot \var Z = \U \cdot \Ni_{c(V)} \cdot \A \cdot \A$.
\end{proof}

\subsection{The first examples of not locally soluble $SI^*$-groups}
The variety we have just constructed is essential for our purposes. Let us use a special notation for it:
$$
	\V^* = \V^*(\U, \V)= \U \cdot \Ni_{c(V)} \cdot \A^2.
$$

\begin{Remark}
Any product of varieties different from the variety $\O$ of all groups also is different from $\O$. Thus, taking into account the fact that we consider only non-trivial word sets $V$, we get that $\V^*(\U, \V)= \O$  if and only if $\U=\O$.
\end{Remark}

Next we need a group $G$ with the following properties: 
\begin{enumerate}
\item[(1)] $G$ is an insoluble $SI^*$-group, 
\item[(2)] $G$ is torsion-free, 
\item[(3)] $G$ is countable, 
\item[(4)] $G$ does not generate the variety of all groups, that is, $G$  satisfies a non-trivial identity. 
\end{enumerate}

The exact form of a group with these properties is not relevant to us. Moreover, later in this section we will be able to build infinitely many groups with properties (1)--(4). What we need now is just {\it one} initial example of such a group. Four different ways of construction of such initial groups with properties (1)--(4) will be given at the end of current section.

The following theorem already provides us with a countable set of $SI^*$-groups, which are not locally soluble.

\begin{Theorem}
\label{Counably_many}
Let $\V$ be an arbitrary variety different from the variety of all groups, that is, $\V$ is defined by a non-trivial word set $V$, and let $\U = \var G$ be a variety generated by a group $G$ satisfying properties (1)--(4). Then there exists a not locally soluble $SI^*$-group $H=H(G,V)$ such that 
$$
	H\in \V^* \quad and \quad H\notin \V,
$$
where $\V^*=\V^*(\U,\V)$ as above.

In particular, it follows from this that there exist countably many not locally soluble $SI^*$-groups, which are pairwise non-isomorphic. 

Moreover, locally soluble $SI^*$-groups mentioned can be constructed to be torsion-free.
\end{Theorem}

\begin{proof}
For any non trivial variety $\V$ let us consider the corresponding (non-trivial) word set $V$ and build the embedding construction of Theorem~\ref{From_That_Paper} for any given group $G$ with properties (1)--(4) and for our $V$. The group $H=H(G,V)$ obtained is an $SI^*$-group. But $H$ cannot be a locally soluble group because it is finitely generated and contains an insoluble subgroup $G$. The fact that the group $H$ belongs to $\V^*$ follows from Corollary~\ref{UNA2}. $H$ does not belong to $\V$ because the verbal subgroup $V(H)$ contains a non trivial subgroup $G$.

Now define $V_0 =V$, $H_0 = G$ and obtain a countable set of examples of not locally soluble $SI^*$-groups $\{ H_i\,|\, i=1,2,\ldots \}$ recursively: 

\begin{enumerate}
\item set the word sets $V_n\sub \Fin$, $i=1,2,\ldots$, to be those corresponding to the varieties 
$$
	\V_n= \V_{n-1} \cup \, \bigl( \var{H_{n-1}} \cdot \Ni_{c(V_{n-1})} \cdot \A^2 \bigr)
$$
(where, clearly, $\V_{n-1}$ 
is the variety corresponding to $(n-1)$th word  set $V_{n-1}$);

\item build the group $H_n=H(H_{n-1}, V_{n-1})$  according to the construction of Theorem~\ref{From_That_Paper}.
\end{enumerate}

It is clear that the groups $H_i$, $i=1,2,\ldots$\,, are not locally soluble two-generator $SI^*$-groups. Moreover, not only $H_i\not\cong H_j$ for any $i\not= j$, but also $\var{H_i}\not= \var{H_j}$, provided that $i\not= j$.

Moreover, if the initial group $H_0=G$ is torsion-free, then all wreath products 
$$
	\bigl((H_{i-1} \Wr S)\Wr C\bigr) \Wr Z, \quad i=1,2,\ldots
$$ 
are torsion-free together with their subgroup $H_i = H(H_{i-1},V_{i-1})$ for $i=1,2,\ldots$
\end{proof}

\subsection{Examples of groups with properties (1)--(4)} 
\label{(1)--(4)}
Now we turn to exact examples of groups $G$ with mentioned properties promised above. 

If some variety $\U\not=\O$ contains a set $\mathcal G =\{G_i \,|\, i=1,2, \ldots ,n,\ldots \}$ of soluble countable groups, such that the solubility lengths of the groups in $\mathcal G$ have no common bound, then the direct product $G$ of elements of $\mathcal G$ clearly satisfies conditions (1), (3) and (4). Moreover, if all groups of $\mathcal G$ are torsin-free, then $G$ also has that property. So what we have to look for, are examples of such sets $\mathcal G$.

I am vary much grateful to Professor B.~I.~Plotkin as well as to Professor  L.~G.~Kov\'acs, who kindly provided me with some examples of mentioned type (the first two of the examples below).  These examples are of different types and they show that in Theorem~\ref{Counably_many} we could build our construction taking one of the following varieties as a `start point' $\U$:  
$$
\A \cdot\B_4,\quad
\A \cdot\K_p \,\,(p\ge 5), \quad
\E_3,  \quad
\A \cdot ( \K_{p} \cap \E_n) \quad (n>3,\,\, p=n+2), 
$$
where $\K_p$ is the Kostrikin variety of all locally finite groups of prime exponent $p$, $\E_n$ is the variety of all $n$-Engel groups.  

\begin{Example}[Kov\'acs]
Let $G_n$ be the free group of rank $n$ in the variety $\A \cdot \B_4$, where $\B_4$ is the Burnside variety of groups of exponent dividing $4$. 
Since the corresponding $\B_4$-free groups $F_n(\B_4)$ are finite, they are soluble. Thus $G_n$ are also soluble. By Razmislov's Theorem~\cite{Razmislov_HH} $\B_4$ is not soluble, so there can be no bound on the solubility lengths of these groups $G_n$. Moreover, the groups $G_n$ are torsion-free; this follows from the result of Kov\'acs~\cite{39Varieties}: relatively free groups of any product variety $\X \cdot \Y$ are torsion-free if and only if relatively free groups of variety $\X$ are torsion-free.
\end{Example}

\begin{Example}[Plotkin]
Since the Kostrikin variety $\K_p$ is a locally nilpotent variety, all relatively free groups $F_n(\K_p)$, $n=1,2,\ldots$, are nilpotent. On the other hand  $\K_p$ is not a soluble variety (see~\cite{RazmislovEngel} for the case $p\ge 5$, and~\cite{BachMach5} for the case $p = 5$). So in order to have a torsion-free group $G$  we can take $\mathcal G$ to be the set of relatively free groups in product variety $\A \cdot \K_p$.
\end{Example}

\begin{Example}
All $3$-Engel groups are locally nilpotent~\cite{Heineken3Engel}. So put $G_n=F_n(\E_3)$. There is no bound on the solubility lengths of these groups, for the free group of infinite rank in variety $\K_5 \cap \E_3$ is insoluble~\cite{BachMach5}.
\end{Example}

\begin{Example}
\label{Example4}
Varieties $\E_n$ for $n>3$ can be another source for examples. Intersections  $\K_{p} \cap \E_n $, where $n\ge3$, $p=n+2$, are locally nilpotent, but still not soluble varieties~\cite{RazmislovEngel}. To get torsion-free groups with desired properties take relatively free groups of finite ranks in the variety $\A \cdot ( \K_{p} \cap \E_n  )$.
\end{Example}

\section{A continuum of not locally soluble $SI^*$-groups}
\label{section3}

Let us first observe, that, in spite of the fact that the set of all varieties of groups (that is, the set of all non-trivial word sets $V$) is of cardinality of continuum, the method of Theorem~\ref{Counably_many} can provide us with only countably many examples of not locally soluble $SI^*$-groups.
The point is that $V$ `takes part' in construction of $H(G,V)$ by means of  the corresponding positive integer $c(V)$, and the number of such integers is, clearly, countable.

\subsection{A group set of Ol'shanskii and the bisections (B)}
We will use a special set of groups built by Ol'shanskii in~\cite{Olshanski_On_problem-of_finite_base} and used for construction of a continuum of varieties of groups.
Namely, let $\{K_n \,|\, n\in \N \}$ be a countable set of finite groups with  following properties:
\begin{enumerate}
\item $K_n \in \G_5 \cap \B_{8qr}$, $n=1,2,\ldots$ , where $\G_5$ is the variety of soluble groups of length at most $5$,  where $q,r$ are different primes and where $\B_{8qr}$ is the Burnside variety of groups of exponents dividing $8qr$;

\item $K_n \notin \var{K_1,\ldots, K_{n-1},K_{n+1}, \ldots}$ for an arbitrary $n\in \N $
\end{enumerate}
(see~\cite{Olshanski_On_problem-of_finite_base} for details). Our final construction will be more economical, if we fix one of the groups presented in Subsection~\ref{(1)--(4)} and continue consideration for that concrete case. So set $G$ be the group of Example~\ref{Example4} (see Remark~\ref{Now_we_can_explane_why_we_used})  and choose the primes $q$ and $r$ to be different from the prime $p=\exp G$.

Let as define a {\it bisection} \bb to be the following representation of the set $\N$ of positive integers
$$
	\N = N' \cup N'',\quad  \text{where $N' \cap N'' = \emptyset $\, and  
	$N',\, N'' \not= \emptyset $.} 
	\leqno {\rm(B)}
$$ 
Define for each bisection (B) a group $K_\bb$ to be the direct product of groups $K_n$, $n\in N'$, a variety $\V_\bb$ to be the variety generated by groups $K_n$, $n\in N''= \N \backslash N'$, and a word set $V_\bb$ to be that corresponding to $\V_\bb$. Denote, further,  
\begin{equation}
\label{Gbb}	
G_\bb = G\times \Fin(\A\cdot \var{K_\bb}).
\end{equation}
It is easy to see that $G_\bb$ is an $SI^*$-group. And since, clearly, $V_\bb$ is a non-trivial word set, we are in position to apply Theorem~\ref{From_That_Paper} to subnormally embed $G_\bb$ into an appropriate two-generator $SI^*$-group $$
	H_\bb = H\bigl(G_\bb, V_\bb \bigr).
$$
The set of all the bisections \bb  is of cardinality of continuum and it would be sufficient to find continuum many different bisections which determine `very different' groups $H_\bb$. 

\begin{Theorem}
\label{technical}
There exists a continuum of the bisections \bb for which the corresponding groups $H_\bb$ are pairwise non-isomorphic and, moreover, they generate pairwise different varieties of groups. 
\end{Theorem}

The proof of this theorem occupies the next subsection.

\subsection{The proof of Theorem~\ref{technical}}
\label{large_subsection}
Assume another bisection $\bbb$ is  given:
$$
	\N = {N'}^* \cup {N''}^*\!\!,\quad  \text{where \,${N'}^* \!\cap {N''}^* = \emptyset $ and  
	${N'}^*\!\!,\, {N''}^* \not= \emptyset $.} 
	\leqno {\hbox{$\bbb$}}
$$ 
The bisections \bb and $\bbb$ are {\it  different} if ${N'}^*\not= N'$ (or equivalently: ${N''}^*\not= N''$).
We begin with the following:
\begin{Lemma}
\label{FirstLemma}
Let the bisections \bb and $\bbb$ be different and let  $G_\bb$  and $G_\bbb$ be the groups constructed according to the rule described above. Then 
$$
\var{G_\bb} \not= \var{G_\bbb}.
$$
\end{Lemma}

\begin{proof}
Since the group 
$$
G \, = \!\!\! \prod_{n=1,\ldots,\infty} \!\!\!F_n\Bigl(\A \cdot ( \K_{p} \cap \E_n  )\Bigr)
$$ 
generates the variety $\A \cdot ( \K_{p} \cap \E_n)$, holds:
\begin{equation*}
\begin{split}
	\var{G_\bb}&=\A \cdot  (\K_{p} \cap \E_n)
	\cup \A\cdot \var{K_\bb}\\
	&=\A \cdot \Bigl( (\K_{p} \cap \E_n )\cup  \var{K_\bb}\Bigr).
\end{split}
\end{equation*}
For arbitrary varieties of groups $\X$, $\Y_1$, $\Y_2$, different from the variety of all groups, $\X \cdot \Y_1= \X \cdot \Y_2$ holds if and only if holds $\Y_1 = \Y_2 $~\cite{HannaNeumann}.
Therefore it is sufficient to show that for different bisections \bb and $\bbb$:
$$
(\K_{p} \cap \E_n )\cup  \var{K_\bb}\not=(\K_{p} \cap \E_n )\cup  \var{K_\bbb}.
$$
Varieties $\var{K_\bb}$ and $\var{K_\bbb}$ are different and thus 
\begin{equation}
\label{1}
	\Fin\bigl(\var{K_\bb}\bigr) \not\cong
	\Fin\bigl(\var{K_\bbb}\bigr).
\end{equation} 
Since varieties $\var{K_\bb}$ and $\var{K_\bbb}$ both have trivial intersection with $\var G$, we get according to~\cite[21.33]{HannaNeumann} that 
\begin{equation}
\label{2}
\begin{split}
	\Fin\bigl(\var{G_\bb}\bigr) &= 
	\Fin\bigl(\var{G} \cup \var{K_\bb} \bigr)\\ 
	&=\Fin\bigl(\var{G} \bigr) \times
	\Fin\bigl(\var{K_\bb} \bigr)
\end{split}
\end{equation} 
and 
\begin{equation}
\label{3}
\begin{split}
	\Fin\bigl(\var{G_\bbb}\bigr) &= 
	\Fin\bigl(\var{G} \cup \var{K_\bbb} \bigr)\\ 
	&=\Fin\bigl(\var{G} \bigr) \times
	\Fin\bigl(\var{K_\bbb} \bigr).
\end{split}
\end{equation} 
Now it follows from~(\ref{2}),~(\ref{3}) and~(\ref{1}) that $\var{G_\bb}$ and $\var{G_\bbb}$ are different varieties.
\end{proof}

Consider for each bisection \bb the complete wreath product
$$
	G_\bb \Wr S_\bb,
$$
where, as in the previous section, $S_\bb$ is defined to be the free nilpotent group of certain rank $r_\bb$ and class $c_\bb$ (and in fact $r_\bb = c_\bb$), such that $S_\bb \notin \V_\bb$. 

\begin{Remark}
\label{which_narrows}
Let us make an observation which narrows the set of the bisections we are working with and shortens our proof. The set of numbers $r_\bb = c_\bb$ is at most countable. Thus there exists a continuum of the bisections determining one and the same value of $c_\bb$ (and, thus, of $r_\bb$). We will show that {\it this} class of bisections already gives birth to a continuum of groups $H_\bb$ generating pairwise different varieties. We replace all values of $c_\bb$ by one $c$ and all values of $r_\bb$ by one $r$ and, thus, we replace all groups $S_\bb$ by single one $S$. (As we will see in Subsection~\ref{variety_containing}, one can simply take $c=r=1$ but this is not relevant for our construction).
\end{Remark}

\begin{Lemma}
Let the bisections \bb and $\bbb$ be different and let $T_\bb = T(G_\bb,S)$ and $T_\bbb = T(G_\bbb,S)$ be groups constructed in analogy with the group $T(G,S)$ in the proof of Theorem~\ref{From_That_Paper}. Then $\var{T_\bb} \not= \var{T_\bbb}$.
\end{Lemma}

\begin{proof}
For arbitrary varieties of groups $\X_1$, $\X_2$, $\Y$, different from variety of all groups, $\X_1 \cdot \Y= \X_2 \cdot \Y$ holds if and only if holds $\X_1 = \X_2 $~\cite{HannaNeumann}.

 Thus, to prove this lemma, it is sufficient to show that
\begin{equation}
\label{4}
\begin{split}
	\var{T_\bb}&=\var{G_\bb} \cdot \Ni_c\\
	\var{T_\bbb}&=\var{G_\bbb} \cdot \Ni_c
\end{split}
\end{equation}
and to apply Lemma~\ref{FirstLemma}. Let us prove, say, the first one of relations~(\ref{4}). Clearly
$$
	\var{T_\bb}\sub \var{G_\bb \Wr S} \sub \var{G_\bb} \cdot \Ni_c.
$$
On the other hand $T_\bb$ contains the first copy of the group $G_\bb$ in the base group of wreath product $G_\bb \Wr S$, as well as the active group $S$. Therefore $T_\bb$ contains a subgroup of the complete wreath product $G_\bb \Wr S$, isomorphic to {\it direct\/} wreath product $G_\bb \Wrr S$ (clearly, this direct wreath product is a proper subgroup of $T_\bb$, for $T_\bb$ contains, say, the element $\chi_g$ with infinite support). Now it remains to notice that $G_\bb \Wrr S$ generates the variety $\var{G_\bb} \cdot \Ni_c$ because the group $S$ {\it discriminates} the variety $\Ni_c$~\cite[35.12, 35.13]{HannaNeumann}.
\end{proof}

The next step of our construction in Section~\ref{section2} was the group $D$. So build here for each bisection \bb a group $D_\bb$. 
It turns out that the commutator subgroup of $D_\bb$ has the same equational theory as the group $T_\bb$:

\begin{Lemma}
If $D_\bb$ and  $T_\bb$ are the same as above, then 
$$
\var{D_\bb}=\var{T_\bb} \A.
$$
In particular, if $\bb$ and $\bbb$ are distinct bisections, then $\var{D_\bb} \not= \var{D_\bbb}$.
\end{Lemma}

\begin{proof}
$D_\bb$ evidently contains the first copy of $T_\bb$ in the base group of the cartesian wreath product $T_\bb \Wr C$. That copy together with the cycle $C$ generate the direct wreath product $T_\bb \Wrr C$ which generates $\var{T_\bb} \A$. It remains to notice that $D_\bb \in \var{T_\bb} \A$ and, thus, $\var{D_\bb} \sub \var{T_\bb} \A$.
\end{proof}

Next we take, as in Section~\ref{section2}, another infinite cycle $Z=\langle  z \rangle$ and form for each bisection \bb, that is, for each group $D_\bb$, a corresponding element $\omega_\bb$ in the complete wreath product $D_\bb \Wr Z$. 
Take for each $D_\bb$ the corresponding two-generator group 
$$
H_\bb= \langle \omega_\bb, z \rangle.
$$

\begin{Remark}
Let us note that, since the elements of  $D_\bb$ can be linearly ordered in many different ways, one can get many different elements $\omega_\bb$. So it is assumed that we have fixed one concrete $\omega_\bb$ for each $D_\bb$, that is, we have fixed some linear order over each group $D_\bb$. 
\end{Remark}

\begin{Lemma}
Let $H_\bb$ be the two-generator group determined by the bisection \bb. Then holds: 
$$
\var{H_\bb}=\var{D_\bb} \A.
$$
In particular, if $\bb$ and $\bbb$ are distinct bisections, then $\var{H_\bb} \not= \var{H_\bbb}$.
\end{Lemma}

\newcommand{\Olshanskii}{The current version of this proof is based on an advice we received from Professor Alexander Yu.~Ol'shanskii. This proof is shorter than the proof we gave in the initial version of this paper.}

\begin{proof}\footnote{\Olshanskii}
Let us take any non-identity $w = w (x_1, \ldots , x_n)$ for the variety $\varr{D_\bb} \A$ and show that $w$ can be falsified on some elements of $H_\bb$. This will prove the point because $H_\bb$ evidently belongs to $\varr{D_\bb} \A$. Take $c_1, \ldots , c_n \in D_\bb \Wrr Z$ such that $w(c_1, \ldots , c_n) \not= 1$. Clearly: $c_i = z^{m_i} \rho_i$, where $\rho_i$ belongs to the base subgroup $D_\bb^{Z}$; $i = 1, \ldots, n$. Finitely many elements $\rho_i$ in this direct wreath product have only finitely many non-trivial ``coordinates'' $\rho_i(z^j)$, $i = 1, \ldots, n$; $j \in \Z$. 
This means that there is a big enough positive integer $n^*$ such that if we replace (trivial) ``coordinates'' $\rho_i(z^j)$, $|j|>n^*$ of each $\rho_i$ by arbitrarily chosen values from the group $D_\bb$, and denote these new strings by $\rho_i'$ correspondingly, then we will still have:
$$
w (z^{m_1} \rho_1', \ldots , z^{m_n} \rho_n',) \not= 1.
$$
Since all the powers $z^{m_1}, \ldots , z^{m_n}$ already belong to $H_\bb$, the proof will be completed if we show that for arbitrary positive integer $n_0$ and arbitrary pregiven values 
$l_{j}\in D_\bb$, $j=-n_0, \ldots , n_0$ the group $H_\bb$ contains such an element $\rho'' \in H_\bb \cap D_\bb^Z$ for which: 
$\rho''(z^{j})=l_{j}$; $j=-n_0, \ldots , n_0$. 

Taking into account the ``shifting'' effect of the element $z$, it will be sufficient to show that for any pregiven $l\in D_\bb$ there is an element $\rho_l''' \in H_\bb \cap D_\bb^{Z}$ such that: 
$$
\rho_l'''(z^j)= \begin{cases}
	l & \text{if $j = 0$} \\
	1 & \text{if $-2 n_0 \le j \le 2 n_0$ and $j \not= 0$}.
	\end{cases}
$$
The elements $\rho''$ will then be products of elements of type $\rho_l'''$ (for various $l$'s and of their conjugates by powers of $z$).
We have: 
$$
\omega_\bb (z^i) = 
\left\{
                  \begin{array}{ll}
                  d_k,  & {\rm if}\,\,\,  i = 2^k, \,\,\, k=0,1,2,...,  \\
                  1,   &  {\rm if}\,\,\,  i \in \Z \backslash \{ 2^k\,\, | \,\,\,  k=0,1,2,...  \}.
                  \end{array}
                  \right.
$$
where this countable group $D_\bb$ is presented as $D_\bb = \{d_0,d_1, \ldots \}$. The element $l$ can be presented as a product $d_i \cdot d_j$ for infinitely many pairs $d_i , d_j \in D_\bb$. On the other hand, the number of all possible pairs $d_i , d_j$ with a common upper bound on $|i|$ and $|j|$ clearly is finite. Thus, there necessarily exists a pair $d_i , d_j$ such that $l = d_i \cdot d_j$ and $2^i, 2^j > 2n_0$. Then:
$$
	\rho'''_l = \omega^{z^{-2^i}} \omega^{z^{-2^j}}.
$$\end{proof}

\subsection{A variety containing all groups $H_\bb$}
\label{variety_containing}

Finally let us find a variety containing all groups $H_\bb$ constructed. For any \bb holds:
\begin{equation*}
\begin{split}
	G_\bb &\in  \bigl(\A \cdot ( \K_{p} \cap \E_n  )\bigr) \cup \bigl(\A \cdot (\G_5 \cap \B_{8qr})\bigr)\\
  	&=\A \cdot \bigl(( \K_{p} \cap \E_n  ) \cup (\G_5 \cap \B_{8qr})\bigr)=\W
\end{split}
\end{equation*}

Further, the group $T_\bb$ belongs to the variety $\W \cdot \Ni_c$, the group $D_\bb$ belongs to the variety $\W \cdot \Ni_c \cdot \A$ and the group $H_\bb$ belongs to the variety $\W \cdot \Ni_c \cdot \A \cdot \A$ for each bisection \bb\!.

\begin{Remark}
Turning back to Remark~\ref{which_narrows}, we note that every one of our varieties $\V_\bb$ is a subvariety of $\G_5 \cap \B_{8qr}$ and, thus, each $V_\bb$ has a consequence of form $x^{8qr}$. So, if $S=F_1(\Ni_1)$ is an infinite cycle, then we already have $\var{S}=\Ni_1 \not\sub \V_\bb$, and we can take $c=1$. So for any \bb the group $T_\bb$ belongs to  $\W \cdot \A$, the group $D_\bb$ belongs to
$
\W \cdot \A \cdot \A =\W \cdot \A^2
$,
and the two-generator group $H_\bb$ lies in $\W \cdot \A^3$. 
\end{Remark}

\begin{Lemma}
For any bisection \bb the two-generator group $H_\bb$ constructed above belongs to the variety
\begin{equation}
\label{containing}
\A \cdot \bigl(( \K_{p} \cap \E_n  ) \cup (\G_5 \cap \B_{8qr})\bigr)\cdot \A^3,
\end{equation}
where $n\ge3$, $p=n+2$ is a prime, and where $q,r$ are primes different from the prime $p$.
\end{Lemma}

\begin{Remark}
\label{Now_we_can_explane_why_we_used}
Now we can explain, why in our construction we prefered the Example~\ref{Example4}. Firstly, the fact that $\exp G$ is coprime to $\exp{(\G_5 \cap \B_{8qr})}$ slightly shortens our original proof (the same is true concerning the group of Example~\ref{Example4}, that is, the group $G=\prod_{i=1}^\infty F_n(\K_p)$), and, secondly, we get a smaller variety for (\ref{containing}) because, $\K_{p} \cap \E_n \sub \K_{p}$.
\end{Remark}

\subsection{A continuum of not locally soluble $SI^*$-groups}
An immediate consequence of Theorem~\ref{technical} and of consideration of  the  previous subsection is the following:

\begin{Theorem}
\label{central}
There exists a continuum of torsion free, not locally soluble two-generator $SI^*$-groups, which generate pairwise different varieties of groups. Moreover, the groups mentioned can be chosen to belong to the variety (\ref{containing}). 
\end{Theorem}

Taking into account the fact the set of all varieties is of cardinality of continuum, we get:

\begin{Theorem}
\label{many_varieties}
The set of all varieties of groups generated by torsion free, not locally soluble two-generator $SI^*$-groups is of the same cardinality as the set of all varieties of groups. Moreover, continuum many varieties of mentioned type are contained in the variety (\ref{containing}). 
\end{Theorem}

If we drop the requirement for the groups $H_\bb$ to be {\it torsion free}, we can slightly change our construction in order to obtain a `smaller value' for variety (\ref{containing}). Namely, take $G$ to be not the same as in Example~\ref{Example4}, but to be the direct product of free groups of finite ranks in variety $ \K_{p} \cap \E_n$:
$$
	G_0=\prod_{n=1,\ldots,\infty}\!\! F_n( \K_{p} \cap \E_n);
$$
take further  $G_\bb$ be not as in (\ref{Gbb}), but:
$$
	G_\bb = G_0 \times \Fin(\var{K_\bb}).
$$
Then a continuum of different bisections \bb still determine a continuum of different insoluble $SI^*$-groups $G_\bb$, which can be used in Theorem~\ref{From_That_Paper} in order to obtain a continuum of groups $H_\bb$. The proof of Lemma~\ref{FirstLemma} can be easily modified (in fact, shortened) for this situation. The rest of the proof in Subsection~\ref{large_subsection} remains unchanged. Since in this case we begin with a group $G_\bb$ not in the variety $\A \cdot ( \K_{p} \cap \E_n  ) \cup \A \cdot (\G_5 \cap \B_{8qr})$, but in the variety $( \K_{p} \cap \E_n  ) \cup  (\G_5 \cap \B_{8qr})$, the groups $H_\bb$ will belong to the variety
\begin{equation}
\label{new_containing}
\bigl(( \K_{p} \cap \E_n  ) \cup (\G_5 \cap \B_{8qr})\bigr)\cdot \A^3
\end{equation}
So we get the following analogs of Theorem~\ref{central} and of Theorem~\ref{many_varieties} respectively:

\begin{Theorem}
\label{new_central}
There exists a continuum of not locally soluble two-gene\-rator $SI^*$-groups, which generate pairwise different varieties of groups. Moreover, the groups mentioned can be chosen to belong to the variety 
(\ref{new_containing}). 
\end{Theorem}

And:

\begin{Theorem}
\label{new_many_varieties}
The set of all varieties of groups generated by not locally soluble two-gene\-rator $SI^*$-groups is of the same cardinality as the set of all varieties of groups. Moreover, continuum many varieties of mentioned type are contained in the variety (\ref{new_containing}). 
\end{Theorem}

\section{Some related questions}
\label{section4}

\subsection{Examples for some other classes of generalized soluble groups}
\label{Examples_for_some_other_classes_of_generalized}
The main tool that enabled us to build examples of not locally soluble $SI^*$-groups is the construction of verbal embedding of countable $SI^*$-groups into two-generator $SI^*$-groups adopted from~\cite{OnNormalEmbeddingsOfSoluble} and described in Section~\ref{section2}. The results of~\cite{OnNormalEmbeddingsOfSoluble} however concern not only $SI^*$-groups but also other classes of generalized soluble groups:
\begin{equation}
\label{12}
	\text{$SI$-, \,$SN$-, \,$SN^*$-, 
	\,$SI$\`{\hskip0.1mm}- \,and \,$SN$\`{\hskip0.1mm}-{\hskip0.1mm}groups.}
\end{equation}

Since there are different ways of definitions of classes of generalized soluble groups (and even different notations of these groups), let us briefly recall the stuff.

A soluble subnormal system of $G$ is the set 
\begin{equation}
\label{system}
	\{G_\delta;\,\, \delta \in \Delta\}
\end{equation}
of subgroups of $G$ with following properties: (i) it contains \1 and $G$, (ii) is linearly ordered by inclusion, (iii) is {\it closed} (that is, contains unions and intersections of its elements), (iv) satisfies the condition $G_\delta \normal G_\delta^\wedge$, where $G_\delta^\wedge$ is the intersection of all elements of~(\ref{system}) greater than $G_\delta$, and (v) every factor $G_\delta^\wedge/G_\delta$ is abelian.
If all subgroups $G_\delta$ are normal in $G$,~(\ref{system}) is said to be a soluble normal system of $G$. The group $G$ is a $SN$-group ($SI$-group) if it posses a soluble subnormal (normal) system.
If moreover the subgroups  $G_\delta$, $\delta \in \Delta$, are well-ordered by inclusion $\sub$, then~(\ref{system}) is a subnormal (normal) soluble ascending series, and $G$ is an $SN^*$-group (an $SI^*$-group) of some ordinal length $\alpha$. If, in the contrary, the subgroups  $G_\delta$, $\delta \in \Delta$, are well-ordered by $\supseteq$, then~(\ref{system}) is a subnormal (normal) soluble descending series, and $G$ is an \,$SN$\`{\hskip0.1mm}-{\hskip0.1mm}group (or a \,$SI$\`{\hskip0.1mm}-{\hskip0.1mm}group) of length $\beta$, where $\beta$ is some inverse ordinal~\cite{RobinsonFinitenessConditions,KuroshThoreyOfGroups,Skornyakov}.

Any countable group $G$ from each one of the classes~(\ref{12}) for each non-trivial word set $V$  is (subnormal and `verbally') embeddable into a two-generator group $H=H(G,V)$ from the very same class~\cite{OnNormalEmbeddingsOfSoluble}. This together with the following lemma give a method of construction of `many' groups from one of the classes~(\ref{12}) that do {\it not} belong to a smaller one of that class (say, examples of $SI$-groups which are not $SI^*$-groups, etc.).

\begin{Lemma}[\cite{RobinsonFinitenessConditions,KuroshThoreyOfGroups}]
Each one of the classes~(\ref{12}) is closed under the operation of taking subgroups.
\end{Lemma}

We simply take, say, one $SI$-group which is not an $SI^*$-group $G$ (e. g. the absolutely free group $G=F$ of finite or countable rank~\cite{RobinsonFinitenessConditions,KuroshThoreyOfGroups}) and for any non-trivial word set $V$ embed it into the corresponding two-generator group $H=H(G,V)$. The latter is still an $SI$-group but not an $SI^*$-group for $G=F \sub H$. The same scheme works for other classes~(\ref{12}). Moreover, the groups constructed can also be not locally soluble if we take an insoluble initial group $G$.

We restrict ourselves to this outline of possible constructions firstly because in the literature examples of mentioned nature are not as rare as that mentioned at the beginning of this paper, and secondly because such a consideration would demand a dozen of pages with proofs very similar to those above varying only in details.

\subsection{Examples for linearly ordered groups}
Linearly ordered (or {\it fully} ordered in terminology of, say,~\cite{EmbeddingsIntoOrdered}) generalized soluble groups have been considered repeatedly (see for instance~\cite{On_non-strictly_simple} and references given there). Thus it is of some interest to note that the groups of Section~\ref{section2} and of examples outlined above can be linearly ordered. 

If the initial group $G$ is linearly ordered, its order can be lifted to that of the group $T(G,S)$ using the fact that the free nilpotent group $S$ can also be linearly ordered. Then the order obtained is being lifted to linear orders of $D(T,C)$ and of $H(G,V)$ using the natural order of cyclic groups $C$ and $Z$. See~\cite{SubnormalEmbedOfOrderedGr,OnNormalEmbeddingsOfSoluble} for details.

The following basic idea is essential to that proofs. Our construction is based on cartesian wreath products and according to the result of Neumann~\cite{EmbeddingsIntoOrdered} a cartesian wreath product of two (non-trivial) groups can never be linearly ordered. This is one of reasons explaining the choice of elements, say, $\chi_g$: on each step not the whole wreath product ($G \Wr S$, $T \Wr C$ or $D \Wr Z$) but a `small part' of it is being used.

\end{document}